\documentclass[10pt,leqno]{amsart}
\usepackage{indentfirst, amssymb,amsmath,amsthm}
\usepackage{csquotes}
\usepackage{graphicx}
\newcommand{\be}{\begin{equation}}
\newcommand{\ee}{\end{equation}}
\newcommand{\beas}{\begin{eqnarray*}}
\newcommand{\eeas}{\end{eqnarray*}}
\newcommand{\bea}{\begin{eqnarray}}
\newcommand{\eea}{\end{eqnarray}}

\numberwithin{equation}{section}
\begin{document}
\title[A Visual Proof that $\pi^{e}<e^{\pi}$]{A Visual Proof that $\pi^{e}<e^{\pi}$}
\date{}
\author[B. Chakraborty ]{Bikash Chakraborty }
\date{}
\address{Department of Mathematics, Ramakrishna Mission Vivekananda
Centenary College, Rahara, West Bengal 700 118, India. }
\email{bikashchakraborty.math@yahoo.com, bikash@rkmvccrahara.org}
\maketitle
\footnotetext{2010 Mathematics Subject Classification: Primary 00A05, Secondary 00A66.}
\textbf{Claim:}  $\pi^{e}<e^{\pi}$.\\
\begin{figure}[h]
\includegraphics[scale=.45]{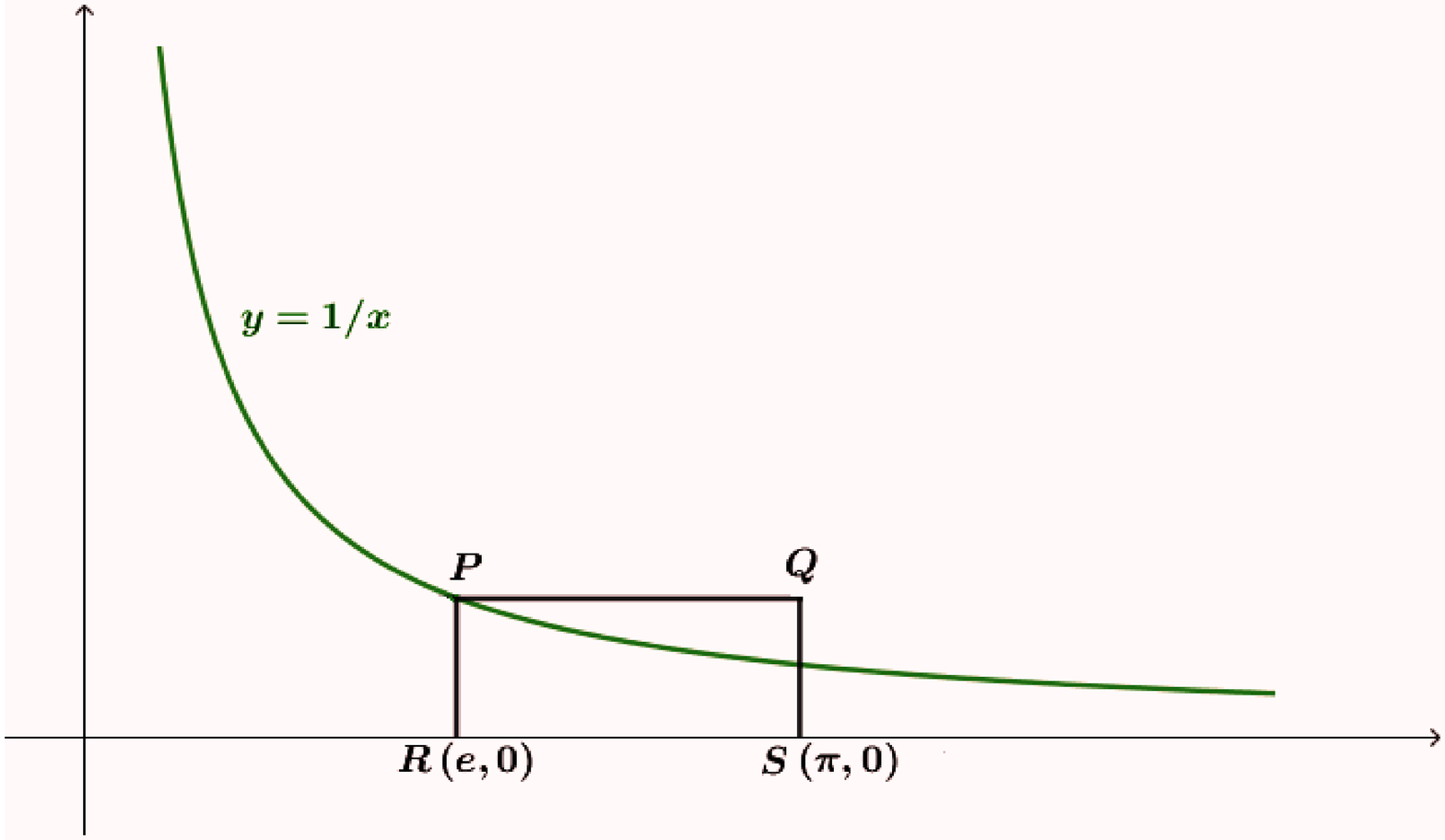}
\centering
\end{figure}
\beas  \ln \pi -1=\int_{e}^{\pi}\frac{dx}{x} &<& \frac{1}{e}(\pi-e)=\frac{\pi}{e}-1,\\
\pi^{e}&<&e^{\pi}.\eeas

\end{document}